\documentclass[MSNbibl,number,citesort,dvips]{arxbj}
\usepackage{upgreek}
\usepackage{graphicx}
\aid{0}
\volume{19}
\issue{5A}
\pubyear{2013}
\firstpage{2010}
\lastpage{2032}
\doi{10.3150/12-BEJ441} 

\makeatletter
\def\pi{\uppi}
\newcommand{\angler}{\rangle}
\newcommand{\anglel}{\langle}
\newcommand{\rrVert}{\Vert}
\newcommand{\rrvert}{\vert}
\newcommand{\llVert}{\Vert}
\newcommand{\llvert}{\vert}
\newtheorem{theorem}{Theorem}[section]
\newremark{example}[theorem]{Example}
\newproclaim{definition}[theorem]{Definition}
\newremark{remark}[theorem]{Remark}
\newtheorem{proposition}[theorem]{Proposition}
\newtheorem{lemma}[theorem]{Lemma}

\newcommand{\bbc}{\mathbb{C}}
\newcommand{\bbr}{\mathbb{R}}
\newcommand{\bbp}{\mathbb{P}}
\newcommand{\bbe}{\mathbb{E}}
\newcommand{\bbf}{\mathbb{F}}
\newcommand{\bbn}{\mathbb{N}}
\newcommand{\bbt}{\mathbb{T}}

\newcommand{\cf}{\mathcal{F}}

\newcommand{\pfe}{\longmapsto}
\newcommand{\pff}{\longrightarrow}
\renewcommand{\backslash}{\setminus}
\renewcommand{\pi}{\uppi}
\newcommand{\eqref}[1]{(\ref{#1})}
\makeatother

\begin{document}
\begin{frontmatter}

\title{Generalization of the Blumenthal--Getoor index to the class of
homogeneous diffusions with jumps and some applications}
\runtitle{Generalized Blumenthal--Getoor index}

\begin{aug}
\author{\fnms{Alexander} \snm{Schnurr}\corref{}\ead[label=e1]{alexander.schnurr@math.tu-dortmund.de}}%
\runauthor{A. Schnurr} 
\address{TU Dortmund, Faculty of Mathematics, Vogelpothsweg 87, 44227
Dortmund, Germany.\\\printead{e1}}
\end{aug}

\received{\smonth{2} \syear{2012}.}

%
\begin{abstract}
We introduce the probabilistic symbol for the class of homogeneous
diffusions with jumps (in the sense of Jacod/Shiryaev).
This concept generalizes the well-known characteristic exponent of a L\'
evy process. Using the symbol, we introduce eight
indices which generalize the Blumenthal--Getoor index $\beta$ and the
Pruitt index $\delta$.
These indices are used afterwards to obtain growth and H\"older
conditions of the process. In the future,
the technical main results will be used to derive further fine
properties. Since virtually all examples of homogeneous
diffusions in the literature are Markovian, we construct a process
which does not have this property.
\end{abstract}

%
\begin{keyword}
\kwd{COGARCH process}
\kwd{Feller process}
\kwd{fine continuity}
\kwd{fine properties}
\kwd{generalized indices}
\kwd{It\^o process}
\kwd{semimartingale}
\kwd{symbol}
\end{keyword}

\end{frontmatter}

\section{Introduction}

Two of the main tools in order to analyze and describe L\'evy processes
are the characteristic exponent and the Blumenthal--Getoor index. In the
present paper, we show that there exist analogous of these concepts for
a much wider class of processes, namely homogeneous diffusions with
jumps (h.d.w.j.) in the sense of Jacod and Shiryaev (\cite{jacodshir},
Definition III.2.18). These indices are used to derive growth and H\"
older conditions for the paths of the process.

A L\'evy process $X$ is a stochastic process with stationary and
independent increments which has a.s. c\`adl\`ag paths (cf. \cite
{sato}). It is a well-known fact that the characteristic function of
$X_t$ can be written as
%
\begin{equation}
\label{charexp} \varphi_{X_t}(\xi)=\bbe^0 \mathrm{e}^{\mathrm{i} X_t'\xi}
= \mathrm{e}^{-t\psi(\xi)},
\end{equation}
where the characteristic exponent $\psi\dvtx \bbr^d\to\bbc$ is a
continuous negative definite function (c.n.d.f.) in the sense of
Schoenberg (cf. \cite{bergforst}, Chapter 2). In fact, one obtains by
the relation \eqref{charexp} a one-to-one correspondence between the
class of c.n.d.f.'s and L\'evy processes. The Blumenthal--Getoor index
was first introduced in \cite{blumenthalget61} in order to analyze H\"
older conditions, the $\gamma$-variation and the Hausdorff-dimension
of the paths of L\'evy processes.

The idea of the present paper is to use the state-space dependent right
derivative at $t=0$ of the characteristic function to obtain the symbol
$p$ of the process which generalizes the characteristic exponent of a
L\'evy process.
The formula reads as follows (for details, see Definition \ref
{defsymbol} below): for $x,\xi\in\bbr^d$
\[
p(x,\xi):=- \lim_{t\downarrow0}\frac{\bbe^x \mathrm{e}^{\mathrm{i}(X^\sigma
_t-x)'\xi}-1}{t},
\]
where $\sigma$ is the first-exit time of a compact neighborhood of $x$.
Since for every fixed $t > 0$, the function $\xi\mapsto\bbe^x
\mathrm{e}^{\mathrm{i}(X_t^\sigma-x)'\xi}$ is the characteristic function of the random
variable $X_t^\sigma-x$ it is continuous and positive definite. By
Corollary 3.6.10 of \cite{niels1}, we conclude that $\xi\pfe-(\bbe^x \mathrm{e}^{\mathrm{i}(X_t^\sigma-x)'\xi} -1)$ is a continuous negative definite
function. Dividing by $t$ preserves this property since the c.n.d.f.'s
form a convex cone. By Lemma 3.6.7 of \cite{niels1}, the above limit
is a negative definite function which is continuous if the convergence
is locally uniform. The idea to analyze objects of this type was
proposed first in \cite{Jaco1998} in the context of universal Markov processes.

We have thus shown that the symbol is a state-space dependent c.n.d.f.
Therefore, we can define and analyze eight indices along the same lines
as in Schilling's article \cite{schilling98} where the case of rich
Feller processes was analyzed. These are Feller processes with the
property that the test functions $C_c^\infty(\bbr^d)$ are contained
in the domain of their generator. The multiplier in the Fourier
representation of the generator of such a process is also a state-space
dependent c.n.d.f. (cf. Example \ref{exfeller} below and for details
the monograph by Jacob \cite{niels1,niels2,niels3}). For these
c.n.d.f.'s, we write $q(x,\xi)$ to distinguish them from the $p(x,\xi
)$ above. In order to introduce and use the indices, Schilling needed
the following two conditions \eqref{growth} and \eqref{sector} which
we state here since they play a role in our considerations, too.
The growth condition is fulfilled, if there exists a $c>0$ such that
\renewcommand{\theequation}{G}
\begin{equation}
\label{growth} \bigl\llVert q(\cdot,\xi) \bigr\rrVert_\infty\leq c
\bigl(1+\llVert \xi \rrVert^2\bigr)
\end{equation}
for every $\xi\in\bbr^d$. The sector condition, which is needed only
for some of the results, is fulfilled, if there exists a $c_0>0$ such
that for every $x,\xi\in\bbr^d$
\renewcommand{\theequation}{S}
\begin{equation}
\label{sector} \bigl\llvert \Im\bigl(q(x,\xi)\bigr) \bigr\rrvert \leq
c_0 \Re\bigl(p(x,\xi)\bigr).
\end{equation}

In \cite{mydiss}, we have shown that every rich Feller process is an
It\^o process in the sense of Cinlar, Jacod, Protter and Sharpe (cf.
\cite{vierleute}, Section 7), that is, a Hunt semimartingale with
characteristics of the form
\renewcommand{\theequation}{\arabic{equation}}
\setcounter{equation}{1}
%
\begin{eqnarray}
\label{hdwj} B_t^{(j)}(\omega) &=&\int_0^t
\ell^{(j)}\bigl(X_s(\omega)\bigr) \, \mathrm{d}s,\qquad
j=1,\ldots,d,
\nonumber
\\
C_t^{jk}(\omega) &=&\int_0^t
Q^{jk}\bigl(X_s(\omega)\bigr) \, \mathrm{d}s, \qquad j,k=1,\ldots,d,
\\
\nu(\omega;\mathrm{d}s,\mathrm{d}y) &=&N\bigl(X_s(\omega),\mathrm{d}y\bigr) \, \mathrm{d}s,
\nonumber
\end{eqnarray}
where for every $x\in\bbr^d$ $\ell(x)$ is a vector in $\bbr^d$,
$Q(x)$ is a positive semi-definite matrix and $N$ is a Borel transition
kernel such that $N(x,\{0\})=0$.
The triplet $(\ell(x),Q(x), N(x,\mathrm{d}y))$ appears in the symbol again (cf.
Theorem \ref{thmstoppedsymbol}). Since the characteristics describe the
local dynamics of the process, it is not surprising that the symbol, as
well as the associated indices, contain a lot of information about the
global and the path properties of the process, like conservativeness
(cf. \cite{schilling98pos}, Theorem 5.5), strong $\gamma$-variation
(cf. \cite{sdesymbol}, Corollary 5.10) or Hausdorff-dimension (cf.
\cite{schilling98hdd}, Theorem 4). By now, all results of this type
were restricted to rich Feller processes. The above considerations show
that It\^o processes would be a natural candidate to generalize the
results on symbols, indices and fine properties. In the present paper
we go even one step further: semimartingales having characteristics of
the form \eqref{hdwj} are called h.d.w.j. It is this class we are
dealing with. In Section \ref{sec2}, we have included an example of this kind,
which is not a Markov process. Philosophically speaking we show that
the symbol, as well as the derived indices, are a concept related to
the underlying semimartingale structure rather than the property of
being memoryless. To this end, new techniques of proof had to be developed.



Here and in the following, we mean by a stochastic process a family of
processes $(X,\bbp^x)_{x\in\bbr^d}$ which is normal, that is, $\bbp^x(X_0=x)=1$. Such a process is called a martingale, continuous, \ldots \,iff
it is w.r.t. every $\bbp^x$ $(x\in\bbr^d)$ a martingale,
continuous, $\ldots\,$. A~stochastic basis $(\Omega,\bbf=(\cf_t)_{t\geq
0},\cf,\bbp^x)_{x\in\bbr^d}$ is always meant to be in the
background. We assume that the usual hypotheses are satisfied.

Before closing this section, we give an overview on what was known
before the present paper. We consider the following classes of processes:
%
\begin{equation}
\label{display} %
\matrix{ %
\mbox{symmetric}
\cr
\alpha
\mbox{-stable} } %
\quad \subseteq\quad \mbox{L\'evy} \quad
\subseteq\quad %
\matrix{ \mbox{rich}
\cr
\mbox{Feller} } %
\quad \subseteq \quad\mbox{It\^o} \quad\subseteq \quad\mbox{h.d.w.j.}
\end{equation}
The symbol was generalized to It\^o processes in \cite{mydiss}. The
indices were known for rich Feller processes satisfying \eqref{growth}
and \eqref{sector}. Fine properties were obtained for the same class,
sometimes under additional assumptions (cf. \cite{schilling98hdd}).
Let us mention that even in the known case of rich Feller processes we
generalize Schilling's results: instead of \eqref{growth} we only need
a local version of this property which is automatically fulfilled by
every rich Feller process.

Let us give a brief outline on how the paper is organized: in the
subsequent section we show that there exists a h.d.w.j. which is not
Markovian. In particular the last inclusion in \eqref{display} is
strict. In Section \ref{sec3} we present the definitions and main results.
Complementary results and several examples, including the COGARCH
process which is used to model financial data, are contained in Section
\ref{sec4}. The proofs are postponed to Section \ref{sec5}, since they are rather
technical. Our main results are Theorems \ref{thmstoppedsymbol}, \ref
{thmttoinfty} and \ref{thmtto0}.

The notation we are using is more or less standard. Vectors are column
vectors. Transposed vectors or matrices are denoted by $'$. Vector
entries are written as follows: $v=(v^{(1)},\ldots,v^{(d)})'$. In the
context of semimartingales we follow mainly \cite{jacodshir}.
Multivariate stochastic integrals are always meant componentwise. This
is true for integrals w.r.t. processes as well as for those w.r.t.
random measures. A function $\chi\dvtx \bbr^d\to\bbr$ is called cut-off
function if it is Borel measurable, with compact support and equal to
one in a neighborhood of zero. In this case $h(y):= \chi(y)\cdot y$ is
a truncation function in the sense of \cite{jacodshir}. Finally, let
$\bbn:=\{0,1,\ldots\}$.

\section{A non-Markovian homogeneous diffusion}\label{sec2}

Virtually all examples of homogeneous diffusions (with or without
jumps) in the literature are Markov processes. Here we construct an
example which is not Markovian.

\begin{example}
We use the construction principle for deterministic processes which we
introduced in \cite{detmp1} and generalized in \cite{detmp2}. Let
$\bbt$ denote the unit sphere in $\bbr^2$.

Within the set $((0,1)'+\bbt)\cup((0,-1)'+\bbt)$, we consider the
following ODE on $[0,\infty[$:
\begin{eqnarray*}
y_1' &=&1-y_2,\qquad y_2'=y_1\qquad
\mbox{for } y_2\geq0,
\\
y_1' &=&y_2+1,\qquad y_2'=-y_1\qquad
\mbox{for } y_2<0
\end{eqnarray*}
with the initial value $y(0)=(y_1(0),y_2(0))'=(0,0)'$ having the
(non-unique) solution
\[
y(t)=\sum_{n\in2\bbn}\pmatrix{\sin(t)
\cr
1-\cos(t)} \cdot
1_{[2n\pi,2(n+1)\pi[}(t) + \pmatrix{\sin(t)
\cr
\cos(t)-1} \cdot1_{[2(n+1)\pi,2(n+2)\pi[}(t).
\]
For the readers convenience, we include the following picture:
%
\begin{figure}[h]

\includegraphics{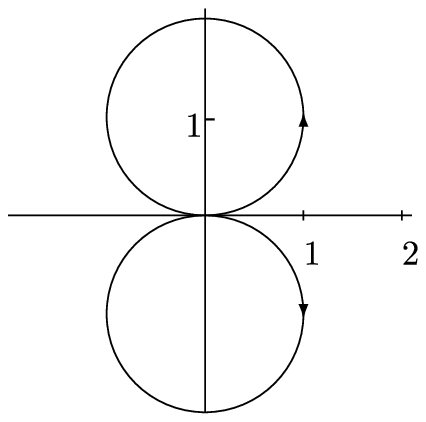}

\end{figure}

We denote by $\widetilde{y}$ the restriction of $y$ to $[0,4\pi[$. On
this interval the function is bijective. The process $X$ is defined as
follows: under the law $\bbp^x$ we have
\[
X_t:= %
\cases{ y\bigl(\widetilde{y}^{-1}(x)+t
\bigr), \quad & for $x\in \bigl((0,1)'+\bbt\bigr)\cup\bigl((0,-1)'+
\bbt\bigr)$,
\cr
x,\quad & else.} %
\]
This process is not Markovian, since
\[
\bbp^{(0, 2)'} \biggl(X_{2\pi}=\pmatrix{0
\cr
-2}\Big\rrvert
X_{\pi
}=\pmatrix{0
\cr
0} \biggr)=1\neq0 =\bbp^{(0,-2)'}
\biggl(X_{2\pi}=\pmatrix{0
\cr
-2}\Big\rrvert X_{\pi
}=\pmatrix{0
\cr
0} \biggr).
\]
On the other hand, $X$ is a homogeneous diffusion with $\ell$ given by
\[
\ell(x)= %
\cases{ \pmatrix{1-x^{(2)}
\cr
x^{(1)}}, \quad &
if $x\in \bigl((0,1)'+\bbt\bigr)$,
\cr
\pmatrix{x^{(2)}+1
\cr
x^{(1)}}, & if $x\in\bigl((0,-1)'+\bbt\bigr)
\backslash{\bigl\{(0,0)'\bigr\}}$,
\cr
0, & else.} %
\]
Anticipating an important concept of the next section, let us mention
that $\ell$ is not continuous on $\bbr^2$, but it is $X$-finely
continuous (cf. Definition \ref{deffinelycont}).
\end{example}

\section{Definitions and main results}\label{sec3}

We have decided to postpone the proofs to Section \ref{sec5}.

\begin{definition}
A \emph{homogeneous diffusion with jumps} (\emph{h.d.w.j.}, for
short) $(X,\bbp^x)_{x\in\bbr^d}$ is a semimartingale with
characteristics of the form
%
\begin{eqnarray}
B_t^{(j)}(\omega) &=& \int_0^t
\ell^{(j)}\bigl(X_s(\omega)\bigr) \, \mathrm{d}s,\qquad j=1,\ldots,d,
\nonumber
\\
C_t^{jk}(\omega) &=& \int_0^t
Q^{jk}\bigl(X_s(\omega)\bigr) \, \mathrm{d}s,\qquad j,k=1,\ldots,d,
\\
\nu(\omega;\mathrm{d}s,\mathrm{d}y) &=& N\bigl(X_s(\omega),\mathrm{d}y\bigr) \, \mathrm{d}s
\nonumber
\end{eqnarray}
for every $x\in\bbr^d$ with respect to a fixed cut-off function $\chi
$. Here $\ell(x)=(\ell^{(1)}(x),\ldots,\ell^{(d)}(x))'$ is a vector in
$\bbr^d$, $Q(x)$ is a positive semi-definite matrix and $N$ is a Borel
transition kernel such that $N(x,\{0\})=0$. We call $\ell$, $Q$ and
$n:=\int_{y\neq0} (1\wedge\llVert  y \rrVert^2) N(\cdot,\mathrm{d}y)$ the
\emph
{differential characteristics} of the process.
\end{definition}

\begin{remark} In the monograph \cite{jacodshir}, this class of
processes is called homogeneous diffusion with jumps, but even there
this name was qualified as `misleading', since the term `diffusion' is
often used for continuous Markov processes: a diffusion with jumps is
not continuous and in Section \ref{sec2} we have seen that it does not have to
be Markovian. However, we decided to stick to the classical name, since
it has become canonical.
\end{remark}

In our considerations, it turned out that the most general assumption
on the differential characteristics, under which we are able to prove
our main results, read as follows.

\begin{definition} \label{deffinelycont}Let $X$ be a h.d.w.j. and
$f\dvtx \bbr^d\to\bbr$ be a Borel-measurable function. $f$ is called
$X$\emph{-finely continuous} (or \emph{finely continuous}, for short)
if the function
%
\begin{equation}
\label{rightcont} t\mapsto f(X_t)=f\circ X_t
\end{equation}
is right continuous at zero $\bbp^x$-a.s. for every $x\in\bbr^d$.
\end{definition}

\begin{remark}
\begin{longlist}[(b)]
\item[(a)] In the context of Markov processes, fine continuity is introduced
differently (see \cite{blumenthalget}, Section II.4, and \cite
{fuglede}). By Theorem 4.8 of \cite{blumenthalget}, this is equivalent
to \eqref{rightcont}.
\item[(b)] If the differential characteristics are continuous, the condition
stated in Definition \ref{deffinelycont} is obviously fulfilled,
since the paths of $X$ are c\`adl\`ag.
\end{longlist}
\end{remark}

The other important assumption on the differential characteristics is
that they are locally bounded. By Lemma 3.3 of \cite
{symbolkillingopenset}, this is equivalent to the local version of the
growth condition: for every compact set $K\subseteq\bbr^d$ there
exists a constant $c_K>0$ such that
\renewcommand{\theequation}{LG}
%
\begin{equation}
\label{localgrowth} \bigl\llvert p(x,\xi) \bigr\rrvert \leq c_K
\bigl(1+\llVert \xi \rrVert^2\bigr)
\end{equation}
for every $x\in K$. This condition is fulfilled by every rich Feller
process (Lemma 3.3 of~\cite{symbolkillingopenset}).

\renewcommand{\theequation}{\arabic{equation}}
\setcounter{equation}{5}

\begin{definition} \label{defsymbol}
Let $X$ be a h.d.w.j., which is conservative and normal, that is, $\bbp^x(X_0=x)=1$. Fix a starting point $x$ and define $\sigma=\sigma^x_k$
to be the first exit time from a compact neighborhood $K:=K_x$ of $x$:
\[
\sigma:=\inf\bigl\{t\geq0 \dvt X_t^x \notin K \bigr\}.
\]

For $\xi\in\bbr^d$, we call $p\dvtx \bbr^d\times\bbr^d\to\bbc$ given by
%
%
\begin{equation}
\label{stoppedsymbol} p(x,\xi):=- \lim_{t\downarrow0}\bbe^x
\frac{\mathrm{e}^{\mathrm{i}(X^\sigma
_t-x)'\xi}-1}{t}
\end{equation}
the \emph{symbol} of the process, if the limit exists and coincides
for every choice of $K$.
\end{definition}

In Example \ref{exfeller}, we show that this symbol coincides with
the classical functional analytic symbol in the case of rich Feller
process. This motivates the name.


\begin{theorem} \label{thmstoppedsymbol}
Let $X$ be a h.d.w.j. such that the differential characteristics $\ell
$, $Q$ and $n$ are locally bounded and finely continuous.
In this case, the limit \eqref{stoppedsymbol} exists and the symbol of
$X$ is
%
\begin{equation}
\label{symbol} p(x,\xi)=-\mathrm{i}\ell(x)'\xi+ \frac{1}{2}
\xi'Q(x) \xi-\int_{y\neq0} \bigl(\mathrm{e}^{\mathrm{i}y'\xi}-1
-\mathrm{i}y'\xi\cdot\chi(y) \bigr) N(x,\mathrm{d}y).
\end{equation}
\end{theorem}


\begin{remark}
\begin{longlist}[(b)]

\item[(a)] If the differential characteristics are continuous, the conditions
of the theorem are fulfilled.

\item[(b)] If the differential characteristics are globally bounded, that is,
if \eqref{growth} is satisfied, the limit \eqref{stoppedsymbol}
without stopping time exists and coincides with the above limit (the
proof is similar).

\item[(c)] Let us mention that the symbol of a L\'evy process is just its
characteristic exponent, that is, $p(x,\cdot)=\psi(\cdot)$ for every
$x\in\bbr^d$. Further examples can be found in the next section.
\end{longlist}
\end{remark}

Now, we define the following helpful quantities for $x\in\bbr^d$ and $R>0$:
%
\begin{eqnarray}
H(x,R)&:=& \sup_{\llVert  y-x \rrVert \leq2R} \sup_{\llVert
\varepsilon \rrVert \leq1} \biggl\llvert p \biggl(y,
\frac{\varepsilon}{R} \biggr) \biggr\rrvert,
\\
H(R)&:=& \sup_{y\in\bbr^d} \sup_{\llVert  \varepsilon \rrVert \leq
1} \biggl\llvert p \biggl(y,
\frac{\varepsilon}{R} \biggr) \biggr\rrvert,
\\
\label{h} h(x,R)&:=& \inf_{\llVert  y-x \rrVert \leq2R} \sup_{\llVert
\varepsilon \rrVert \leq1} \Re p \biggl(y,
\frac{\varepsilon}{4\kappa R} \biggr),
\\
\label{h2} h(R)&:=& \inf_{y\in\bbr^d} \sup_{\llVert  \varepsilon \rrVert \leq1} \Re p \biggl(y,
\frac{\varepsilon}{4\kappa R} \biggr).
\end{eqnarray}
In \eqref{h} and \eqref{h2} $\kappa=(4 \arctan(1/2 c_0))^{-1}$
where $c_0$ comes from the sector condition \eqref{sector} as defined
in the introduction. In particular, $h(x,R)$ and $h(R)$ are only
defined if \eqref{sector} is satisfied and only in this case they will
be used below.

\begin{definition}
The quantities (cf. \cite{schilling98}, Definitions 4.2 and 4.5)
\begin{eqnarray*}
\beta_0&:=&\sup \Bigl\{\lambda\geq0 \dvt \mathop{\lim\sup}_{R\to\infty}
R^\lambda H(R) =0 \Bigr\},
\\
\underline{\beta_0}&:=&\sup \Bigl\{\lambda\geq0 \dvt
\mathop{\lim\inf}_{R\to
\infty} R^\lambda H(R) =0 \Bigr\},
\\
\overline{\delta_0}&:=&\sup \Bigl\{\lambda\geq0 \dvt
\mathop{\lim\sup}_{R\to
\infty} R^\lambda h(R) =0 \Bigr\},
\\
\delta_0&:=&\sup \Bigl\{\lambda\geq0 \dvt \mathop{\lim\inf}_{R\to\infty}
R^\lambda h(R) =0 \Bigr\}
\end{eqnarray*}
are called \emph{indices of $X$ at the origin}, while
\begin{eqnarray*}
\beta_\infty^x&:=&\inf \Bigl\{\lambda> 0 \dvt
\mathop{\lim\sup}_{R\to0} R^\lambda H(x,R) =0 \Bigr\},
\\
\underline{\beta_\infty^x}&:=&\inf \Bigl\{\lambda> 0 \dvt
\mathop{\lim\inf}_{R\to0} R^\lambda H(x,R) =0 \Bigr\},
\\
\overline{\delta_\infty^x}&:=&\inf \Bigl\{\lambda> 0 \dvt
\mathop{\lim\sup}_{R\to0} R^\lambda h(x,R) =0 \Bigr\},
\\
\delta_\infty^x&:=&\inf \Bigl\{\lambda> 0 \dvt
\mathop{\lim\inf}_{R\to0} R^\lambda h(x,R) =0 \Bigr\}
\end{eqnarray*}
are the \emph{indices of $X$ at infinity}.
\end{definition}

\begin{example} In the case of symmetric $\alpha$-stable processes,
all indices coincide and they are equal to $\alpha$. For so called
stable-like Feller processes (cf. \cite{bass88a,negoro94}) with
uniformly bounded exponential function, that is, $0<\alpha_0 \leq
\alpha(x)\leq\alpha_\infty<1$ one obtains $\beta_0=\underline
{\beta_0}=\alpha_0$ and $\delta_0=\overline{\delta_0}=\alpha_\infty$ (see \cite{schilling98}, Example 5.5). For more examples,
consult the next section.
\end{example}

The following proposition is the key ingredient for using the symbol to
analyze fine properties of a stochastic process. Similar results were
proved for L\'evy processes by Pruitt in \cite{pruitt81} and for rich
Feller processes satisfying \eqref{growth} and \eqref{sector} by
Schilling in \cite{schilling98}. We write
\[
(X_\cdot-x)_t^*:=\sup_{s\leq t} \llVert
X_s-x \rrVert
\]
for the maximum process.

\begin{proposition} \label{proptechmain}
Let $X$ be a h.d.w.j. such that the differential characteristics of $X$
are locally bounded and finely continuous. In this case, we have
%
\begin{equation}
\label{firstestimate} \bbp^x \bigl((X_\cdot-
x)_t^* \geq R \bigr) \leq c_d \cdot t \cdot H(x,R)
\end{equation}
for $t\geq0$, $R>0$ and a constant $c_d>0$ which can be written down
explicitly and only depends on the dimension $d$.

If \emph{\eqref{sector}} holds in addition, we have
%
\begin{equation}
\label{secondestimate} \bbp^x \bigl((X_\cdot-x)_t^*
<R \bigr) \leq c_\kappa\cdot\frac{1}{t} \cdot\frac{1}{h(x,R)}
\end{equation}
for a constant $c_\kappa$ only depending on the $c_0$ of the sector condition.
\end{proposition}


Using this result and standard Borel--Cantelli techniques, we obtain the
following two theorems which describe the behavior of the process at
infinity respective zero.

\begin{theorem} \label{thmttoinfty}
Let $X$ be a h.d.w.j. such that the differential characteristics of $X$
are locally bounded and finely continuous. Then we have
%
\begin{eqnarray}
\lim_{t\to\infty} t^{-1/\lambda}(X_\cdot-x)_t^*&=&0\qquad
\mbox{for all } \lambda< \beta_0,
\\
\mathop{\lim\inf}_{t\to\infty} t^{-1/\lambda}(X_\cdot-x)_t^*&=&0\qquad
 \mbox{for all } \beta_0 \leq\lambda< \underline{
\beta_0}.
\end{eqnarray}
If the symbol $p$ of the process $X$ satisfies \emph{\eqref{sector}}, then we
have in addition
%
\begin{eqnarray}
\mathop{\lim\sup}_{t\to\infty} t^{-1/\lambda}(X_\cdot-x)_t^*&=&
\infty\qquad \mbox {for all } \overline{\delta_0} < \lambda\leq
\delta_0,
\\
\lim_{t\to\infty} t^{-1/\lambda}(X_\cdot-x)_t^*&=&
\infty \qquad\mbox{for all } \delta_0 < \lambda.
\end{eqnarray}
All these limits are meant $\bbp^x$-a.s. with respect to every $x\in
\bbr^d$.
\end{theorem}

\begin{theorem} \label{thmtto0}
Let $X$ be a h.d.w.j. such that the differential characteristics of $X$
are locally bounded and finely continuous. Then we have
%
\begin{eqnarray}
\label{proved} \lim_{t\to0 } t^{-1/\lambda}(X_\cdot-x)_t^*&=&0 \qquad
\mbox{for all } \lambda> \beta_\infty^x,
\\
\mathop{\lim\inf}_{t\to0 } t^{-1/\lambda}(X_\cdot-x)_t^*&=&0\qquad
 \mbox{for all } \beta_\infty^x \geq\lambda> \underline{
\beta_\infty^x}.
\end{eqnarray}
If the symbol $p$ of the process $X$ satisfies \emph{\eqref{sector}}, then
we have in addition
%
\begin{eqnarray}
\mathop{\lim\sup}_{t\to0 } t^{-1/\lambda}(X_\cdot-x)_t^*&=&
\infty\qquad \mbox{for all } \overline{\delta_\infty^x} > \lambda
\geq\delta_\infty^x,
\\
\lim_{t\to0 } t^{-1/\lambda}(X_\cdot-x)_t^*&=&
\infty \qquad \mbox{for all } \delta_\infty^x > \lambda.
\end{eqnarray}
All these limits are meant $\bbp^x$-a.s with respect to every $x\in
\bbr^d$.
\end{theorem}

The relation between indices of this type associated with L\'evy
processes and the classical Blumenthal--Getoor respective Pruitt indices
were analyzed in Section 5 of \cite{schilling98}.

\section{Examples, applications, complementary results}\label{sec4}

In the present section, we show how the above results can be used for
some classes of processes. The first example explains the connection
with the classical Markovian theory. The second one deals with L\'evy
driven SDEs having unbounded coefficients and the third one with the
COGARCH process.

\begin{example} \label{exfeller}
Let $X$ be a Feller processes, that is, a strong Markov process such
that
\begin{longlist}[(F2)]
\item[(F1)] $T_t\dvtx C_\infty(\bbr^d) \to C_\infty(\bbr^d)$ for every $t\geq0$,

\item[(F2)] $\lim_{t\downarrow0} \llVert  T_tu-u \rrVert_\infty=0$ for
every $u\in
C_\infty(\bbr^d)$,
\end{longlist}
where
\[
T_t u(x):= \bbe^x u(X_t),\qquad t\geq0, x\in
\bbr^d
\]
and $C_\infty(\bbr^d)$ denotes the real-valued continuous functions
vanishing at infinity. The generator $(A,D(A))$ of the process is the
closed operator given by
%
\begin{equation}
\label{generator} Au:=\lim_{t \downarrow0} \frac{T_t u -u}{t} \qquad\mbox{for } u
\in D(A),
\end{equation}
where the domain $D(A)$ consists of all $u\in C_\infty(\bbr^d)$ for
which the limit \eqref{generator} exists uniformly. Using a classical
result due to Courr\`ege \cite{courrege}, Jacob (cf. \cite{niels1},
Section 4.5) showed that the generator $A$ of a process of this kind
can be written in the following way:
\[
Au(x)= -\int_{\bbr^d} \mathrm{e}^{\mathrm{i}x'\xi} q(x,\xi) \widehat{u}(
\xi)\, \mathrm{d}\xi \qquad\mbox{for } u\in C_c^\infty\bigl(\bbr^d
\bigr),
\]
where $\widehat{u}(\xi)=(2\pi)^{-d}\int \mathrm{e}^{-\mathrm{i}y'\xi}u(y) \, \mathrm{d}y$ denotes
the Fourier transform. The functional analytic symbol $q\dvtx \bbr^d \times
\bbr^d \to\bbc$ has the following properties:
it is locally bounded,
$q(\cdot,\xi)$ is measurable for every $\xi\in\bbr^d$ and
$q(x,\cdot)$ is a c.n.d.f. for every $x\in\bbr^d$.
The last point means that the symbol admits a `state-space dependent'
L\'evy--Khinchine formula like \eqref{symbol}. In Lemma 3.3 of \cite
{symbolkillingopenset}, we have shown that the symbol $q$ always
satisfies \eqref{localgrowth}.


By Theorem 3.10 of \cite{mydiss}, every rich Feller process is an It\^
o process and the differential characteristics are equal to the L\'evy
triplet of the symbol. From Corollary 4.5 of the same thesis, we deduce
that for a rich Feller process with finely continuous differential
characteristics the functional analytic symbol and the probabilistic
symbol do coincide, that is, $p(x,\xi)=q(x,\xi)$ for every $x,\xi\in
\bbr^d$. Furthermore, this shows that the case treated in Schilling
\cite{schilling98} is encompassed by our considerations. Having a look
at his Theorem 3.5, this does not seem to be the case, because the
characteristics look differently, but this is due to a different choice
of the cut-off function.
\end{example}

\begin{example}
Let $(Z_t)_{t\geq0}$ be an $\bbr^n$-valued L\'evy process. The
solution of the stochastic differential equation
%
\begin{eqnarray}\label{levysde}
\mathrm{d}X_t^x&=&\Phi\bigl(X_{t-}^x
\bigr) \, \mathrm{d}Z_t,
\nonumber
\\
[-8pt]
\\[-8pt]
X_0^x&=&x, \qquad x\in\bbr^d,
\nonumber
\end{eqnarray}
where $\Phi\dvtx  \bbr^d \to\bbr^{d \times n}$ is locally Lipschitz
continuous and satisfies the standard linear growth condition, admits
the symbol
\[
p(x,\xi)=\psi\bigl(\Phi(x)'\xi\bigr),
\]
where $\psi\dvtx \bbr^n\to\bbc$ denotes the characteristic exponent of
the L\'evy process. This was shown in~\cite{sdesymbol}. Fine
properties could only be obtained for the case of bounded $\Phi$,
because in general the solution of the above SDE is not rich Feller.
Using the classical characterization of It\^o processes due to Cinlar
and Jacod (\cite{cinlarjacod81}, Theorem 3.33), it is straightforward
to show that $X$ belongs to this class. Since $\Phi$ and $\psi$ are
continuous, the symbol is finely continuous. Along the same lines as in
\cite{sdesymbol}, we obtain the following two results.
\end{example}

\begin{theorem}
Let $p(x,\xi)$ be a state-space dependent c.n.d.f. which can be
written as $p(x,\xi)=\psi(\Phi(x)'\xi)$ where $\psi\dvtx \bbr^n\to
\bbc$ is a c.n.d.f. and $\Phi\dvtx \bbr^d\to\bbr^{d\times n}$ is
locally Lipschitz continuous and satisfies the linear growth condition.
In this case there exists a corresponding It\^o process, that is, a
process $X$ with symbol $p(x,\xi)$.
\end{theorem}

\begin{theorem}
Let $Z$ be a driving L\'evy process with non-constant symbol. Let $X$
be the solution of \eqref{levysde} such that $d=n$ and the rank of
$\Phi$ is equal to $d$ in every point. Then
\[
\lim_{t\to0} t^{-1/\lambda} (X_\cdot- x)_t^*=0\qquad
\mbox{if } \lambda > \beta_\infty,
\]
where $\beta_\infty$ is the index of the driving L\'evy process $Z$.
\end{theorem}

\begin{example}
Let us recall how the COGARCH process is defined (cf. \cite{cogarch}):

Let $Z=(Z_t)_t$ be a L\'evy process with triplet $(\ell,Q,N)$ and fix
$0<\delta< 1, \beta>0, \lambda\geq0$. The volatility process
$(S_t)_{t\geq0}$ is the solution of the SDE
\begin{eqnarray*}
\mathrm{d}S^2_t &=& \beta\, \mathrm{d}t + S_t^2
\biggl( \log\delta \,\mathrm{d}t + \frac{\lambda
}{\delta} \,\mathrm{d} \biggl(\sum
_{0<s\leq t} (\Delta Z_s)^2 \biggr) \biggr),
\\
S_0 &=& S\, (>0).
\end{eqnarray*}
The process
\[
G_t:=g+\int_0^t S_{s-}
\,\mathrm{d}Z_t,\qquad g\in\bbr
\]
is called \emph{COGARCH process}. The pair $(G_t, S_t)$ is a (normal)
Markov process which is is homogeneous in space in the first component.
It is not a Feller process, at least not a $C_\infty$-Feller process.
Furthermore $(G_t, S_t^2)$ is an It\^o process, which follows by
combining Theorem~3.33 of \cite{cinlarjacod81} with Proposition
IX.5.2 of \cite{jacodshir}. To avoid problems which might arise for
processes defined on $\bbr\times\bbr_+$, we consider the logarithmic
squared volatility, that is, the process $(G_t,V_t)=(G_t, \log(S_t^2))$.
This process admits the symbol $p\dvtx \bbr^2 \times\bbr^2 \pff\bbc$
given by
\begin{eqnarray*}
&&p \biggl( \pmatrix{g\cr
v} ,\xi \biggr)\\
&&\quad = -\mathrm{i} \xi_1 \biggl(\ell \mathrm{e}^{v/2} + \mathrm{e}^{v/2} \int
_{\bbr
\backslash\{ 0\}} y \cdot(1_{\{\llvert  \mathrm{e}^{v/2}y \rrvert <1\}} \cdot 1_{\{ \llvert  \log(1+(\lambda/\delta) y^2) \rrvert <1 \}} -
1_{\{
\llvert  y \rrvert  < 1\}}) N(\mathrm{d}y) \biggr)
\\
&&\qquad {}- \mathrm{i}\xi_2 \biggl(\frac{\beta}{\mathrm{e}^{v}} + \log\delta+ \int
_{\bbr\backslash\{ 0\}} \log\biggl(1+\frac{\lambda}{\delta}y^2\biggr)
\cdot(1_{\{\llvert  \mathrm{e}^{v/2}y \rrvert <1\}} \cdot1_{\{ \llvert  \log
(1+(\lambda /\delta) y^2) \rrvert <1 \}} ) N(\mathrm{d}y) \biggr)
\\
&&\qquad {}+\frac{1}{2} \xi_1^2 \mathrm{e}^{v} Q
\\
&&\qquad{}-\int_{\bbr^2 \backslash\{0 \}} \bigl( \mathrm{e}^{\mathrm{i}(z_1,z_2)\xi}-1-\mathrm{i} z'\xi
\cdot(1_{\{\llvert  z_1 \rrvert <1\}} \cdot1_{\{
\llvert  z_2 \rrvert <1 \}} ) \bigr) \tilde{N} \biggl( \pmatrix{g
\cr
v} , \mathrm{d}z \biggr),
\end{eqnarray*}
where $\tilde{N}$ is the image measure
\[
\tilde{N} \biggl( \pmatrix{g
\cr
v} , \mathrm{d}z \biggr) = N(f_v\in \mathrm{d}z)
\]
under $f\dvtx \bbr\to\bbr^2$ given by
\[
f_v(w)= \pmatrix{\mathrm{e}^{v/2}w
\cr
\log\bigl(1+(\lambda/\delta)
w^2\bigr)}.
\]
This was shown in \cite{cogarchsymbol}. A typical driving term in
mathematical finance is the variance gamma process (cf. \cite
{vargamma} and \cite{cogarch2}). This is a pure jump L\'evy process with
\[
N(\mathrm{d}y)= \frac{C}{\llvert  y \rrvert } \exp \bigl( -(2C)^{-1/2} \llvert y \rrvert
\bigr) \, \mathrm{d}y
\]
for a constant $C>0$. In order to have a concrete example, let $\lambda
=2$, $\delta=1/2$, $\beta=10$ and $C=2$. Using standard calculus we
obtain that $\beta_0=1$. The calculations are elementary but tedious.
By Theorem \ref{thmttoinfty}, we obtain for $g\in\bbr$
\[
\lim_{t\to\infty} t^{-1/\lambda}(G_\cdot-g)_t^*=0\qquad
\mbox{for all } \lambda< 1.
\]
In the future, the indices will be used in order to obtain other fine
properties of non-Feller processes.
\end{example}

Now we consider the special case of a process which consists of
independent components.

\begin{proposition}
Let $X$ be a $d$-dimensional vector of independent  h.d.w.j.'s
$X^{(j)}$ with symbols $p^{(j)}, j=1,\ldots,d$. The process $X$ admits the symbol
\[
p(x,\xi)=p^{(1)}\bigl(x^{(1)},\xi^{(1)}
\bigr)+\cdots+p^{(d)}\bigl(x^{(d)},\xi^{(d)}\bigr).
\]
\end{proposition}

\begin{pf}
We give the proof for two components. The general case follows
inductively. Let $X$ and $Y$ be independent h.d.w.j.'s with symbols
$p(x,\xi_1)$, respectively, $q(y,\xi_2)$, where the sum of the dimensions of
$x$ and $y$ is $d$, and consider: 
\begin{eqnarray*}
&&\bbe^{(x,y)} \frac{ \mathrm{e}^{\mathrm{i}(X_t-x)'\xi_1 + \mathrm{i}(Y_t-y)'\xi_2}-1}{t}
\\
&&\quad= \frac{\bbe^{(x,y)}  (\mathrm{e}^{\mathrm{i}(X_t-x)'\xi_1 + \mathrm{i}(Y_t-y)'\xi_2}
)-1}{t}
\\
&&\quad= \frac{\bbe^{x}  (\mathrm{e}^{\mathrm{i}(X_t-x)'\xi_1} ) \cdot\bbe^{y}
(\mathrm{e}^{\mathrm{i}(Y_t-y)'\xi_2} )-1}{t}
\\
&&\quad= \frac{\bbe^{x}  (\mathrm{e}^{\mathrm{i}(X_t-x)'\xi_1} ) \cdot\bbe^{y}
(\mathrm{e}^{\mathrm{i}(Y_t-y)'\xi_2} ) -\bbe^{y} (\mathrm{e}^{\mathrm{i}(Y_t-y)'\xi_2}
)+\bbe^{y} (\mathrm{e}^{\mathrm{i}(Y_t-y)'\xi_2} ) -1}{t}
\\
&&\quad= \frac{\bbe^{x}  (\mathrm{e}^{\mathrm{i}(X_t-x)'\xi_1} ) -1}{t}\cdot\bbe^{y} \bigl(\mathrm{e}^{\mathrm{i}(Y_t-y)'\xi_2} \bigr) +
\frac{\bbe^{y}
(\mathrm{e}^{\mathrm{i}(Y_t-y)'\xi_2} ) -1}{t}.
\end{eqnarray*}
The three terms on the right-hand side tend to $-p(x,\xi_1)$, $1$ and
$-q(y,\xi_2)$, respectively. Hence, the result.
\end{pf}

\section{Proofs of the main results}\label{sec5}

In this section, we present the proofs of the main results.

\begin{pf*}{Proof of Theorem \protect\ref{thmstoppedsymbol}}
Let $x\in\bbr^d$ and let the stopping time defined as in Definition
\ref{defsymbol} where $K$ is an arbitrary compact neighborhood of
$x$. We give the one dimensional proof, since the multidimensional
version works alike; only the notion becomes more involved. First, we
use It\^o's formula under the expectation and obtain
\renewcommand{\theequation}{\Roman{equation}}
\setcounter{equation}{0}
%
\begin{eqnarray}
&&\frac{1}{t} \bbe^x \bigl( \mathrm{e}^{\mathrm{i}(X^\sigma_t-x)\xi} -1 \bigr)
\nonumber
\\
\label{termone}&&\quad= \frac{1}{t} \bbe^x \biggl(\int_{0+}^t
\mathrm{i} \xi \mathrm{e}^{\mathrm{i}(X^\sigma
_{s-}-x)\xi} \,\mathrm{d}X^\sigma_s \biggr)
\\
\label{termtwo}&&\qquad+ \frac{1}{t} \bbe^x \biggl(\frac{1}{2} \int
_{0+}^t-\xi^2 \mathrm{e}^{\mathrm{i}(X^\sigma_{s-}-x)\xi} \, \mathrm{d}
\bigl[X^\sigma,X^\sigma\bigr]_s^c \biggr)
\\
\label{termthree}&&\qquad+ \frac{1}{t} \bbe^x \biggl(\mathrm{e}^{-\mathrm{i}x\xi} \sum
_{0<s\leq t} \bigl(\mathrm{e}^{\mathrm{i}
\xi X_s^\sigma}-\mathrm{e}^{\mathrm{i}\xi X^\sigma_{s-}}-\mathrm{i}\xi
\mathrm{e}^{\mathrm{i}\xi X^\sigma_{s-}} \Delta X^\sigma_s \bigr) \biggr).
\end{eqnarray}
The left-continuous process $X^\sigma_{t-}$ is bounded on $ [ [0,
\sigma]  ]$. Furthermore, we have $(\Delta X)^\sigma= (\Delta
X^\sigma)$ and $X^\sigma$ admits the stopped characteristics
\renewcommand{\theequation}{\arabic{equation}}
\setcounter{equation}{23}
%
\begin{eqnarray}
\label{stoppedchars} B^\sigma_t(\omega)&=&\int
_{0}^{t\wedge\sigma(\omega)} \ell \bigl(X_s(\omega)
\bigr) \, \mathrm{d}s = \int_0^t \ell\bigl(X_s(
\omega)\bigr) 1_{ [ [0,
\sigma ]  ]}(\omega, s) \, \mathrm{d}s,
\nonumber
\\
C_t^\sigma(\omega)&=&\int_0^t
Q\bigl(X_s(\omega)\bigr) 1_{ [ [0, \sigma
]  ]}(\omega, s) \, \mathrm{d}s,
\\
\nu^\sigma(\omega;\mathrm{d}s,\mathrm{d}y)&:=& 1_{ [ [0, \sigma]  ]}(\omega , s)N\bigl(X_s(\omega),\mathrm{d}y\bigr) \, \mathrm{d}s\nonumber %
\end{eqnarray}
with respect to the fixed cut-off function $\chi$. One can now set the
integrand at the right endpoint of the stochastic support to zero, as
we are integrating with respect to Lebesgue measure:
\begin{eqnarray*}
B^\sigma_t(\omega)&=&\int_0^t
\ell\bigl(X_s(\omega)\bigr) 1_{ [ [0 ,
\sigma [  [}(\omega, s) \, \mathrm{d}s,
\\
C_t^\sigma(\omega)&=&\int_0^t
Q\bigl(X_s(\omega)\bigr) 1_{ [ [0 , \sigma
[  [}(\omega, s) \, \mathrm{d}s,
\\
\nu^\sigma(\omega;\mathrm{d}s,\mathrm{d}y)&=&1_{ [ [0 , \sigma[  [}(\omega , s) N
\bigl(X_s(\omega),\mathrm{d}y\bigr) \, \mathrm{d}s.
\end{eqnarray*}
In the first two lines, the integrand is now bounded, because $\ell$
and $Q$ are locally bounded and $\llVert  X^{\sigma}_s(\omega) \rrVert <k$ on
$ [ 0 , \sigma(\omega)  [$ for every $\omega\in\Omega$.
In what follows, we will deal with the terms one-by-one.
To calculate the first term, we use the canonical decomposition of the
semimartingale (see \cite{jacodshir}, Theorem II.2.34) which we write
as follows
%
\begin{eqnarray}\label{candec}
X_t^\sigma&=&X_0 + X_t^{\sigma,c}
+ \int_0^{t \wedge\sigma} \chi (y)y \bigl(
\mu^{X^\sigma}(\cdot;\mathrm{d}s,\mathrm{d}y)-\nu^\sigma(\cdot ;\mathrm{d}s,\mathrm{d}y) \bigr)\nonumber
\\[-8pt]\\[-8pt]
&&{}+\check{X}^\sigma(\chi)+B_t^\sigma(\chi),
\nonumber
\end{eqnarray}
where $\check{X}_t=\sum_{s\leq t}(\Delta X_s (1- \chi(\Delta X_s)))$.
Therefore, term \eqref{termone} can be written as
\begin{eqnarray*}
&&\frac{1}{t} \bbe^x \biggl(\int_{0+}^t
\mathrm{i} \xi \mathrm{e}^{\mathrm{i}(X^\sigma_{s-}-x)
\xi}\,\mathrm{d}\biggl( \underbrace{X_t^{\sigma,c}}_{\mbox{\scriptsize{(IV)}}}
+ \underbrace{\int_0^{t \wedge\sigma} \chi(y)y \bigl(
\mu^{X^\sigma
}(\cdot;\mathrm{d}s,\mathrm{d}y)-\nu^\sigma(\cdot;\mathrm{d}s,\mathrm{d}y)
\bigr)}_{\mbox{\scriptsize{(V)}}}
\\
&&\hspace*{103pt}{}+\underbrace{\check{X}^\sigma(\chi)}_{\mbox{\scriptsize{(VI)}}} + \underbrace
{B_t^\sigma(\chi) }_{\mbox{\scriptsize{(VII)}}} \biggr) \biggr).
\end{eqnarray*}
We use the linearity of the stochastic integral. Our first step is to
prove for term (IV)
\[
\bbe^x\int_{0+}^t \mathrm{i} \xi
\mathrm{e}^{\mathrm{i}(X_{s-}^\sigma-x)\xi} \,\mathrm{d}X_s^{\sigma,c} = 0.
\]
The integral $\mathrm{e}^{\mathrm{i}(X_{t-}^\sigma-x)\xi} \bullet X_t^{\sigma,c}$ is a
local martingale, since $X_t^{\sigma,c}$ is a local martingale. To see
that it is indeed a martingale, we calculate the following:
\begin{eqnarray*}
\bigl[\mathrm{e}^{\mathrm{i}(X^\sigma-x)\xi} \bullet X^{\sigma,c}, \mathrm{e}^{\mathrm{i}(X^\sigma
-x)\xi} \bullet
X^{\sigma,c} \bigr]_t 
&=& \bigl[\mathrm{e}^{\mathrm{i}(X^\sigma-x)\xi} \bullet X^{c}, \mathrm{e}^{\mathrm{i}(X^\sigma-x)\xi}
\bullet X^{c} \bigr]^\sigma_t
\\
&=& \int_0^t
\bigl(\mathrm{e}^{\mathrm{i}(X^\sigma_s-x)\xi}\bigr)^2 1_{ [ [0, \sigma]
]}(s)\, \mathrm{d}
\bigl[X^{c},X^{c}\bigr]_s
\\
&=& \int_0^t \bigl(
\bigl(\mathrm{e}^{\mathrm{i}(X^\sigma_s-x)\xi}\bigr)^2 1_{ [ [0 , \sigma
[  [}(s)Q(X_s)
\bigr) \, \mathrm{d}s,
\end{eqnarray*}
where we have used several well known facts about the square bracket.
The last term is uniformly bounded in $\omega$ and therefore, finite
for every $t\geq0$.
This means that $\mathrm{e}^{\mathrm{i}(X^\sigma_t-x)\xi} \bullet X_t^{\sigma,c}$ is
an $L^2$-martingale which is zero at zero and therefore, its expected
value is constantly zero.

The same is true for the integrand (V). We show that the function
$H_{x,\xi}(\omega,s,y):=\mathrm{e}^{\mathrm{i}(X_{s-}^\sigma-x)\xi}\cdot y \chi(y)$
is in the class $F_p^2$ of Ikeda and Watanabe (see \cite{ikedawat},
Section II.3), that is,
\[
\bbe^x \int_0^t \int
_{y\neq0} \bigl\llvert \mathrm{e}^{\mathrm{i}(X_{s-}^\sigma-x)\xi} \cdot y\chi(y) \bigr
\rrvert^2 \nu^\sigma(\cdot;\mathrm{d}s,\mathrm{d}y) <\infty.
\]
To prove this, we observe
\begin{eqnarray*}
&&\bbe^x \int_0^t \int
_{y\neq0} \bigl\llvert \mathrm{e}^{\mathrm{i}(X_{s-}^\sigma-x)\xi} \bigr
\rrvert^2 \cdot\bigl\llvert y\chi(y) \bigr\rrvert^2
\nu^\sigma(\cdot;\mathrm{d}s,\mathrm{d}y)
\\
&&\quad= \bbe^x \int_0^t \int
_{y\neq0} \bigl\llvert y\chi(y) \bigr\rrvert^2
1_{ [ [0 , \sigma[  [}(\omega, s) N(X_s,\mathrm{d}y) \, \mathrm{d}s.
\end{eqnarray*}
Since we have by hypothesis $\llVert  \int_{y\neq0}(1\wedge y^2)
1_{ [ [0 , \sigma[  [} N(\cdot,\mathrm{d}y) \rrVert_\infty<
\infty$, this expected
value is finite. Therefore, the function $H_{x,\xi}$ is in $F_p^2$ and
we conclude that
\begin{eqnarray*}
&&\int_0^t \mathrm{e}^{\mathrm{i}(X_{s-}^\sigma-x)\xi}\, \mathrm{d} \biggl(\int
_0^{s\wedge\sigma
} \int_{y\neq0} \chi(y)y
\bigl(\mu^{X^\sigma}(\cdot;\mathrm{d}r,\mathrm{d}y)-\nu^\sigma (\cdot;\mathrm{d}r,\mathrm{d}y)\bigr)
\biggr)
\\
&&\quad=\int_0^t\int_{y\neq0}
\bigl(\mathrm{e}^{\mathrm{i}(X_{s-}-x)\xi} \chi (y)y \bigr) \bigl(\mu^{X^\sigma}(\cdot;\mathrm{d}s,\mathrm{d}y)-
\nu^\sigma(\cdot;\mathrm{d}s,\mathrm{d}y)\bigr)
\end{eqnarray*}
is a martingale. The last equality follows from \cite{jacodshir},
Theorem I.1.30.

Now we deal with the second term \eqref{termtwo}. Here we have
\[
\bigl[X^\sigma,X^\sigma\bigr]_t^c=
\bigl[X^c, X^c\bigr]_t^\sigma=C_t^\sigma=
\bigl(Q(X_{t}) \bullet t\bigr)^\sigma = \bigl(Q(X_t)
\cdot1_{ [ [0 , \sigma[  [}(t)\bigr)\bullet t
\]
and therefore,
\[
\frac{1}{2} \int_{0+}^t-\xi^2
\mathrm{e}^{\mathrm{i}(X^\sigma_{s-}-x)\xi} \,\mathrm{d}\bigl[X^\sigma ,X^\sigma\bigr]_s^c
= - \frac{1}{2} \xi^2 \int_{0}^t
\mathrm{e}^{\mathrm{i}(X^\sigma_{s-}-x)\xi} Q(X_{s}) \cdot1_{ [ [0 , \sigma[  [}(t) \, \mathrm{d}s.
\]
Since $Q$ is finely continuous and locally bounded, we obtain by
dominated convergence
\[
-\lim_{t\downarrow0} \frac{1}{2}\xi^2 \frac{1}{t}
\bbe^x \int_0^t \mathrm{e}^{\mathrm{i} (X_{s}-x)\xi}
Q(X_{s}) 1_{ [ [0 , \sigma[  [}(s) \, \mathrm{d}s =-\frac
{1}{2}
\xi^2 Q(x).
\]
For the finite variation part of the first term, that is, (VII), we
obtain analogously
\[
\lim_{t\downarrow0} \mathrm{i}\xi\frac{1}{t} \bbe^x \int
_0^t \mathrm{e}^{\mathrm{i}
(X_{s}-x)\xi} \ell(X_{s})
1_{ [ [0 , \sigma[  [}(s) \, \mathrm{d}s =\mathrm{i}\xi\ell(x).
\]
Now we have to deal with the various jump parts. At first, we write the
sum in \eqref{termthree} as an integral with respect to the jump
measure $\mu^{X^\sigma}$ of the process:
\begin{eqnarray*}
&& \mathrm{e}^{-\mathrm{i}x\xi} \sum_{0<s\leq t} \bigl(\mathrm{e}^{\mathrm{i}X_s\xi}-\mathrm{e}^{\mathrm{i}X_{s-}\xi}-\mathrm{i}
\xi \mathrm{e}^{\mathrm{i}\xi X_{s-}} \Delta X_s \bigr)
\\
&&\quad= \mathrm{e}^{-\mathrm{i}x\xi} \sum_{0<s\leq t}
\bigl(\mathrm{e}^{\mathrm{i} X_{s-}\xi
}\bigl(\mathrm{e}^{\mathrm{i}\xi\Delta X_s}-1 -\mathrm{i}\xi\Delta X_s\bigr)
\bigr)
\\
&&\quad=\int_{]0,t]\times\bbr^d} \bigl( \mathrm{e}^{\mathrm{i} (X_{s-}-x)\xi
}\bigl(\mathrm{e}^{\mathrm{i}\xi y}-1
-\mathrm{i}\xi y\bigr)1_{ \{ y\neq0  \} } \bigr) \mu^{X^\sigma}(\cdot ;\mathrm{d}s,\mathrm{d}y)
\\
&&\quad= \int_{]0,t]\times \{ y\neq0  \} } \bigl(\mathrm{e}^{\mathrm{i}
(X_{s-}-x) \xi}\bigl(\mathrm{e}^{\mathrm{i}\xi y}-1
-\mathrm{i} \xi y\chi(y)\bigr) \bigr) \mu^{X^\sigma
}(\cdot;\mathrm{d}s,\mathrm{d}y)
\\
&&\qquad+ \int_{]0,t]\times \{ y\neq0  \} } \bigl( \mathrm{e}^{\mathrm{i}
(X_{s-}-x) \xi}\bigl( -\mathrm{i}\xi y \cdot
\bigl(1-\chi(y)\bigr)\bigr) \bigr) \mu^{X^\sigma
}(\cdot;\mathrm{d}s,\mathrm{d}y).
\end{eqnarray*}
The last term cancels with the one we left behind from \eqref
{termone}, given by (VI). For the remainder-term, we get:
\begin{eqnarray*}
&& \frac{1}{t} \bbe^x \int_{]0,t]\times \{ y\neq
0  \} }
\bigl(\mathrm{e}^{\mathrm{i} (X_{s-}-x) \xi}\bigl(\mathrm{e}^{\mathrm{i}\xi y}-1 -\mathrm{i} \xi y\chi(y)\bigr) \bigr)
1_{ [ [0 , \sigma[  [}(\cdot, s) \mu^{X^\sigma}(\cdot ;\mathrm{d}s,\mathrm{d}y)
\\
&&\quad= \frac{1}{t} \bbe^x \int_{]0,t]\times \{ y\neq
0  \} }
\bigl(\mathrm{e}^{\mathrm{i} (X_{s-}-x) \xi}\bigl(\mathrm{e}^{\mathrm{i}\xi y}-1 -\mathrm{i} \xi y\chi(y)\bigr) \bigr)
1_{ [ [0 , \sigma[  [}(\cdot, s) \nu^\sigma(\cdot ;\mathrm{d}s,\mathrm{d}y)
\\
&&\quad= \frac{1}{t} \bbe^x \int_{]0,t]\times \{ y\neq
0  \} } \bigl(
\underbrace{\mathrm{e}^{\mathrm{i} (X_{s-}-x) \xi}\bigl(\mathrm{e}^{\mathrm{i}\xi y}-1 -\mathrm{i} \xi y\chi (y)\bigr) \bigr)
1_{ [ [0 , \sigma[  [}(\cdot, s) N(X_s,\mathrm{d}y) }_{:=g(s-,\cdot
)} \mathrm{d}s
\\
&&\quad= \frac{1}{t} \bbe^x \int_{]0,t]\times \{ y\neq
0  \} }
\bigl(\mathrm{e}^{\mathrm{i} (X_{s-}-x) \xi}\bigl(\mathrm{e}^{\mathrm{i}\xi y}-1 -\mathrm{i} \xi y\chi(y)\bigr) \bigr)
1_{ [ [0 , \sigma[  [}(\cdot, s) N(X_s,\mathrm{d}y) \, \mathrm{d}s.
\end{eqnarray*}
Here we have used the fact that it is possible to integrate with
respect to the compensator of a random measure instead of the measure
itself, if the integrand is in $F_p^1$ (see \cite{ikedawat}, Section
II.3). The function $g(s,\omega)$ is measurable and bounded by our
assumption, since $\llvert  \mathrm{e}^{\mathrm{i}\xi y}-1 -\mathrm{i} \xi y\chi(y) \rrvert  \leq \mathit{const}
\cdot(1\wedge\llVert  y \rrVert^2)$. Hence, $g\in F_p^1$.
Again by bounded convergence, we obtain
\begin{eqnarray*}
&&\lim_{t \downarrow0} \frac{1}{t} \bbe^x \int
_{0}^t \mathrm{e}^{\mathrm{i}(X_{s}-x)
\xi} \int
_{y\neq0} \bigl(\mathrm{e}^{\mathrm{i}y \xi}-1-\mathrm{i}y \xi\chi(y) \bigr)
N(X_{s},\mathrm{d}y) \, \mathrm{d}s
\\
&&\quad= \int_{y\neq0} \bigl(\mathrm{e}^{\mathrm{i} y \xi}- 1 - \mathrm{i}y \xi\chi (y)
\bigr) N(x, \mathrm{d}y).
\end{eqnarray*}
This is the last part of the symbol. Here, we have used the continuity
assumption on $N(x,\mathrm{d}y)$.
\end{pf*}


Now we prepare the proof of Proposition \ref{proptechmain}, our
technical main result. It will turn out to be useful to have a closer
look at the symbol \eqref{symbol}. The real part of $p$ is $\Re
(p(x,\xi))= (1/2)\xi'Q(x)\xi-\int_{y\neq0} (\cos(y'\xi)-1)
N(x,\mathrm{d}y)$ and therefore, we obtain
%
\begin{equation}
\label{cos} \int_{y\neq0} \bigl(1-\cos\bigl(y'
\eta\bigr)\bigr) N(x,\mathrm{d}y) \leq\Re\bigl(p(x,\xi)\bigr).
\end{equation}
We assume for the remainder of this section:
$R>0$ and $S>2R$. $\chi$ is a fixed cut-off function such that
\[
\chi\in C_c^\infty\bigl(\bbr^d\bigr) ;\qquad
1_{B_R(0)} \leq\chi\leq1_{B_{2R}(0)} ;\qquad \chi(y)=\chi(-y)\qquad \mbox{for
every } y\in\bbr^d.
\]
The stopping time $\sigma=\sigma_R$ is defined as follows
\[
\sigma:= \inf\bigl\{t\geq0 \dvt  \llVert X_t-x \rrVert > S\bigr\}.
\]
We need the following two lemmas.
%
\begin{lemma}\label{lemft}
For every $z\in\bbr^d$, we have
\[
\bigl(\llVert z \rrVert^2 \wedge1\bigr) \leq c \bigl(1-\mathrm{e}^{-\llVert  z \rrVert ^2/2}
\bigr) \leq c\bigl(\llVert z \rrVert^2 \wedge1\bigr),
\]
where $c=1/(1-\exp(-1/2))$ and
\[
\bigl(1-\mathrm{e}^{-\llVert  z \rrVert ^2/2} \bigr)=\int_{\bbr^d} \bigl(1-\cos
\bigl(z'\eta\bigr)\bigr) h_d \,\mathrm{d}\eta
\]
with
\[
h_d(\eta)=\frac{1}{(\sqrt{2\pi})^d} \mathrm{e}^{-\llVert  \eta \rrVert
^2/2}.
\]
\end{lemma}

The proof is elementary and hence omitted.

\begin{lemma} \label{lemoutsourced}
Let $p(x,\xi)$ be the symbol \eqref{symbol} and $R > 0$. Then we have
\[
\int_{z\neq0} \biggl( \biggl\llVert \frac{z}{2R} \biggr
\rrVert^2 \wedge1 \biggr) N(y,\mathrm{d}z) \leq\widetilde{c}_d
\sup_{\llVert  \varepsilon \rrVert
\leq1} \biggl\llvert p \biggl(y,\frac{\varepsilon}{2R} \biggr) \biggr
\rrvert,
\]
where $\widetilde{c}_d=2c(d+1)$ with the $c$ of Lemma \ref{lemft}.
\end{lemma}

\begin{pf}
By the above lemma, we obtain
\begin{eqnarray*}
\mathit{LHS} &\leq& c\int_{z\neq0} \biggl(1-\exp \biggl( - \biggl\llVert
\frac
{z}{2R} \biggr\rrVert^2 \Big\slash 2 \biggr) \biggr) N(y,\mathrm{d}z)
\\
&=& c\int_{z\neq0} \int_{\bbr^d} \biggl( 1-\cos
\biggl( \frac
{1}{2R} \bigl(z'\eta\bigr) \biggr) \biggr)
\frac{1}{(\sqrt{2\pi})^d} \mathrm{e}^{-\llVert  \eta \rrVert ^2/2}\, \mathrm{d}\eta N(y,\mathrm{d}z)
\\
&\leq& c\int_{\bbr^d} \Re p \biggl(y,\frac{\eta}{2R} \biggr)
h_d(\eta ) \,\mathrm{d}\eta
\\
&\leq&2c \int_{\bbr^d} \sup_{\llVert  \varepsilon \rrVert  \leq1} \biggl\llvert p
\biggl(y,\frac{\varepsilon}{2R} \biggr) \biggr\rrvert \bigl(1+\llVert \eta
\rrVert^2\bigr) h_d(\eta)\, \mathrm{d}\eta
\\
&=& \sup_{\llVert  \varepsilon \rrVert  \leq1} \biggl\llvert p \biggl(y,\frac {\varepsilon}{2R} \biggr)
\biggr\rrvert \int_{\bbr^d} 2c \bigl(1+\llVert \eta
\rrVert^2\bigr) h_d(\eta)\, \mathrm{d}\eta,
\end{eqnarray*}
where we have used the Tonelli--Fubini theorem, the inequality \eqref
{cos} and a standard estimate of the c.n.d.f. $\eta\mapsto p(y,\eta
/(2R))$ as it can be found in the proof of Lemma 3.2 in~\cite{mydiss}.
\end{pf}

\begin{pf*}{Proof of Proposition \ref{proptechmain}}
Let $X$ be a h.d.w.j. such that the differential characteristics $(\ell
,Q,n)$ of $X$ are locally bounded and finely continuous. At first, we
show that for $S,R$ and $\sigma$ as above we have
%
\begin{equation}
\label{techmainstopped} \bbp^x\bigl(\bigl(X_\cdot^\sigma-x
\bigr)_t^* \geq2R\bigr) \leq c_d \cdot t \cdot
\sup_{\llVert  y-x \rrVert \leq S} \sup_{\llVert  \varepsilon \rrVert
\leq1} \biggl\llvert p \biggl(y,
\frac{\varepsilon}{2R} \biggr) \biggr\rrvert,
\end{equation}
where $c_d=4d+16\widetilde{c}_d$. Having proved this the result
follows easily.

The semimartingale characteristics of the stopped process $X^\sigma$
are given in \eqref{stoppedchars} above. Now, we use a double stopping
technique introducing
\[
\tau_R:= \inf\bigl\{t\geq0\dvt  \bigl\llVert \Delta
X_t^\sigma \bigr\rrVert >R \bigr\}.
\]
We start with
%
\begin{equation}
\label{separateone} \bbp^x \bigl( \bigl(X_\cdot^\sigma-x
\bigr)_t^* \geq2R \bigr) \leq\bbp^x \bigl(
\bigl(X_\cdot^\sigma- x\bigr)_t^* \geq2R,
\tau_R > t \bigr) + \bbp^x (\tau_R \leq t )
\end{equation}
and deal with the terms on the right-hand side one after another,
starting with the first one.

We show how to separate the first term of \eqref{separateone} again in
order to get control over the big jumps. Let $\check{X}$ be as defined
in equation \eqref{candec}. The semimartingale $\check{X}^\sigma$
admits the following third characteristic: $\chi(y) 1_{ [ [0,
\sigma ]  ]} (s) N(X_s,\mathrm{d}y) \, \mathrm{d}s$.
Now let $u=(u_1,\ldots,u_d)'\dvtx \bbr^d \to\bbr^d$ be such that $u_j \in
C_b^2(\bbr^d)$ is 1-Lipschitz continuous, $u_j$ depends only on
$x^{(j)}$ and is zero in zero for $j=1,\ldots,d$.
We define the auxiliary process
\[
\check{M}_t:=u\bigl(\check{X}_t^\sigma-x\bigr) -
\int_0^{t\wedge\sigma} F_s \,\mathrm{d}s,
\]
where
%
\begin{eqnarray}\label{auxprocess}
 F_s^{(j)} &=& \partial_j u(
\check{X}_{s-}-x) \ell^{(j)}(X_{s-})
\nonumber
\\
&&{}-\frac{1}{2} \partial_j\, \partial_j u(
\check{X}_{s-}-x) Q^{jj} (X_{s-})
\nonumber
\\[-8pt]
\\[-8pt]
&&{}-\int_{z\neq0} \bigl( u(\check{X}_{s-}-x+z)-u(
\check{X}_{s-}-x)
\nonumber
\\
&& \hspace*{34pt}{}- \chi(z) z^{(j)}\, \partial_j u(\check{X}_{s-}-x)
\bigr) \chi(z) N(X_{s-}, \mathrm{d}z).
\nonumber
\end{eqnarray}
$\check{M}$ is a local martingale by \cite{jacodshir}, Theorem II.2.42
and by Lemma 3.7 of \cite{mydiss} we have under~(\ref{localgrowth}):
\[
\bigl\llvert F_s^{(j)} \bigr\rrvert \leq \mathit{const}\cdot \sum
_{0\leq\llvert  \alpha
\rrvert \leq2} \bigl\llVert \partial^\alpha u \bigr
\rrVert_\infty
\]
since $u_j\in C_b^2(\bbr^d)$. In particular, for every fixed $t>0$
$\check{M}$ is an $L^2$-martingale on $[0,t]$.
Now we define
\[
D:= \biggl\{\omega\in\Omega\dvt \int_0^{t\wedge\sigma(\omega)} \bigl
\llVert F_s(\omega) \bigr\rrVert \,\mathrm{d}s \leq R \biggr\}
\]
and obtain
%
\begin{equation}
\label{separatetwo} \bbp^x \bigl( \bigl(X_\cdot^\sigma-
x\bigr)_t^* \geq2R, \tau_R > t \bigr) \leq
\bbp^x \bigl( \bigl(X_\cdot^\sigma- x
\bigr)_t^* \geq2R, \tau_R > t, D \bigr) +
\bbp^x \bigl(D^c\bigr).
\end{equation}
%
Using Doob's inequality and the Lipschitz property of $u$, we obtain at first
\begin{eqnarray*}
\bbp^x \bigl( u\bigl(X_\cdot^\sigma- x
\bigr)_t^* \geq2R, \tau_R > t, D \bigr) &\leq&
\bbp^x \biggl( u\bigl(X_\cdot^\sigma- x
\bigr)_t^*-\int_0^{\cdot\wedge
\sigma}
F_s \,\mathrm{d}s \geq R, \tau_R > t, D \biggr)
\\
&\leq&\bbp^x\bigl(\check{M}_{t\wedge\sigma}^* \geq R\bigr)
\\
&\leq&\frac{1}{R^2} \bbe^x \bigl(\bigl\llVert
\check{M}_t^\sigma \bigr\rrVert^2\bigr)
\\
&\leq&\frac{1}{R^2} \sum_{j=1}^d
\bbe^x \bigl( \bigl[\check{X}_\cdot^{(j)},
\check{X}_\cdot^{(j)}\bigr]_t^\sigma
\bigr).
\end{eqnarray*}
Since
\[
\bbe^x \bigl( \bigl[\check{X}_\cdot^{(j)},
\check{X}_\cdot^{(j)}\bigr]_t^\sigma \bigr)
= \bbe^x \bigl( \bigl\anglel \check{X}_\cdot^{(j),c},
\check{X}_\cdot^{(j),c}\bigr\angler_t^\sigma
\bigr) + \bbe^x \biggl( \int_0^{t\wedge\sigma}
\int_{z\neq0} \bigl(z^{(j)}\bigr)^2
\chi(z)^2 N(X_s,\mathrm{d}z) \, \mathrm{d}s \biggr)
\]
we obtain
\begin{eqnarray*}
&&\bbp^x \bigl( u\bigl(X_\cdot^\sigma- x
\bigr)_t^* \geq2R, \tau_R > t, D \bigr)
\\
&&\quad \leq\frac{1}{R^2} \sum_{j=1}^d
\bbe^x \int_0^{t\wedge\sigma}
Q^{jj}(X_s) \, \mathrm{d}s + \bbe^x \int
_0^{t\wedge\sigma} \int_{z\neq0}
\frac{\llVert  z
\rrVert^2}{R^2} \chi(z)^2 N(X_s,z) \, \mathrm{d}s
\\
&&\quad \leq4\sum_{j=1}^d \bbe^x
\int_0^{t\wedge\sigma} \biggl(\frac
{e^{\prime}_j}{2R}
Q(X_s) \frac{e_j}{2R} \biggr) \, \mathrm{d}s + 4^2
\bbe^x \int_0^{t\wedge\sigma} \int
_{z\neq0} \biggl( \biggl\llVert \frac{z}{2R} \biggr
\rrVert^2 \wedge1 \biggr) N(X_s,\mathrm{d}z) \, \mathrm{d}s
\\
&&\quad \leq4t \sum_{j=1}^d
\sup_{s<t\wedge\sigma} \Re p \biggl(X_s, \frac{e_j}{2R} \biggr) +
4^2 \sup_{\llVert  y-x \rrVert \leq S} \int_0^{t\wedge\sigma}
\int_{z\neq0} \biggl( \biggl\llVert \frac{z}{2R} \biggr
\rrVert^2 \wedge1 \biggr) N(y,\mathrm{d}z) \, \mathrm{d}s
\\
&&\quad \leq4td \sup_{\llVert  y-x \rrVert  \leq S} \sup_{\llVert
\varepsilon \rrVert  \leq
1} \biggl\llvert p \biggl(y,
\frac{\varepsilon}{2R} \biggr) \biggr\rrvert + 4^2 t\sup_{\llVert  y-x \rrVert \leq S}
\widetilde{c}_d \sup_{\llVert  \varepsilon \rrVert  \leq1} \biggl\llvert p \biggl(y,
\frac{\varepsilon}{2R} \biggr) \biggr\rrvert,
\end{eqnarray*}
where we have used Lemma \ref{lemoutsourced} on the second term. By
choosing a sequence $(u_n)_{n\in\bbn}$ of functions of the type
described above which tends to the identity in a monotonous way, we obtain
%
\begin{equation}
\label{finalterm1D} \bbp^x \bigl( \bigl(X_\cdot^\sigma-
x\bigr)_t^* \geq2R, \tau_R > t, D \bigr) \leq \bigl(4d+
4^2 \widetilde{c}_d \bigr) t \sup_{\llVert  y-x \rrVert  \leq S}
\sup_{\llVert  \varepsilon \rrVert  \leq1} \biggl\llvert p \biggl(y,\frac
{\varepsilon}{2R} \biggr) \biggr
\rrvert .
\end{equation}
%
Now we deal with the second term of \eqref{separatetwo}. By the Markov
inequality, we get
\[
\bbp^x\bigl(D^c\bigr)=\bbp^x \biggl( \int
_0^{t\wedge\sigma} \llVert F_s \rrVert \,\mathrm{d}s > R
\biggr) \leq\frac{1}{R} \sum_{j=1}^d
\bbe^x \biggl(\int_0^{t\wedge\sigma} \bigl
\llvert F_s^{(j)} \bigr\rrvert \,\mathrm{d}s \biggr)=:(*).
\]
Again, we chose a sequence $(u_n)_{n\in\bbn}$ of functions as we
described in \eqref{auxprocess}, but this time it is important that
the first and second derivatives are uniformly bounded. Since the $u_n$
converge to the identity, the first partial derivatives tend to 1 and
the second partial derivatives to 0. In the limit ($n\to\infty$), we obtain
%
\begin{eqnarray}
(*)&\leq&\frac{1}{R} \sum_{j=1}^d
\bbe^x \int_0^{t \wedge\sigma} \biggl\llvert
\ell^{(j)}(X_s) + \int_{z\neq0}
\bigl(-z^{(j)} \chi(z) + \bigl(\chi (z)\bigr)^2
z^{(j)} \bigr) N(X_s,\mathrm{d}z) \biggr\rrvert \,\mathrm{d}s
\nonumber
\\[-2pt]
&\leq&2\sum_{j=1}^d \bbe^x
\int_0^{t\wedge\sigma} \biggl\llvert \frac
{\ell^{(j)}(X_s)}{2R} +
\int_{z\neq0} \sin \biggl(\frac{z'e_j}{2R} \biggr) -
\frac {z^{(j)}\chi(z)}{2R} N(X_s,\mathrm{d}z) \biggr\rrvert \,\mathrm{d}s \label {termA}
\\[-2pt]
&&{}+2\sum_{j=1}^d \bbe^x \int
_0^{t\wedge\sigma} \biggl\llvert \int_{z\neq 0}
\frac{(\chi(z))^2 z^{(j)}}{2R} - \sin \biggl(\frac
{z'e_j}{2R} \biggr) N(X_s,\mathrm{d}z)
\biggr\rrvert \,\mathrm{d}s. \label{termD}
\end{eqnarray}
For term \eqref{termA}, we get
%
\begin{eqnarray}
&&2 \sum_{j=1}^d \bbe^x\int
_0^{t\wedge\sigma} \biggl\llvert \frac{\ell
(X_s)'e_j}{2R} + \int
_{z\neq0} \sin \biggl(\frac{z'e_j}{2R} \biggr) -
\frac{z'e_j \chi(z)}{2R} N(X_s,\mathrm{d}z) \biggr\rrvert \,\mathrm{d}s
\nonumber
\\[-2pt]
&&\quad\leq2td\sup_{s\leq t\wedge\sigma} \bbe^x \biggl\llvert \frac{\ell
(X_s)'e_j}{2R}
+ \int_{z\neq0} \sin \biggl(\frac{z'e_j}{2R} \biggr) -
\frac{z'e_j \chi(z)}{2R} N(X_s,\mathrm{d}z) \biggr\rrvert
\\[-2pt]
&&\quad\leq2td \sup_{\llVert  y-x \rrVert \leq S} \sup_{\llVert
\varepsilon \rrVert \leq1} \biggl\llvert \Im p \biggl(y,
\frac{\varepsilon}{2R} \biggr) \biggr\rrvert \label {finalterm1DcA}\nonumber
\end{eqnarray}
and for term \eqref{termD}
%
\begin{eqnarray}\label {finalterm1DcB}
&&2\sum_{j=1}^d \bbe^x \int
_0^{t\wedge\sigma} \biggl\llvert \int_{z\neq
0}
\frac{(\chi(z))^2 z'e_j}{2R}- \sin \biggl(\frac{z'e_j}{2R} \biggr) N(X_s,\mathrm{d}z)
\biggr\rrvert \,\mathrm{d}s
\nonumber
\\[-2pt]
&&\quad\leq2\sum_{j=1}^d \bbe^x
\int_0^{t\wedge\sigma} \biggl\llvert \int
_{B_{2R}(0)\backslash\{0\}} 1-\cos \biggl(\frac{z'e_j}{2R} \biggr)
N(X_s,\mathrm{d}z) \biggr\rrvert
\nonumber
\\[-9pt]\\[-9pt]
&&\qquad +\biggl\llvert \int_{B_{2R}(0)^c} 1 N(X_s,\mathrm{d}z) \biggr
\rrvert\, \mathrm{d}s
\nonumber
\\[-2pt]
&&\quad\leq2td \sup_{\llVert  y-x \rrVert \leq S} \sup_{\llVert
\varepsilon \rrVert  \leq1} \Re p \biggl(y,
\frac{\varepsilon}{2R} \biggr) + 2^2 td\sup_{\llVert
y-x \rrVert \leq S}
\widetilde{c}_d \sup_{\llVert  \varepsilon
\rrVert  \leq1} \biggl\llvert p \biggl(y,
\frac{\varepsilon}{2R} \biggr) \biggr\rrvert,\nonumber
\end{eqnarray}
where we have used again Lemma \ref{lemoutsourced} on the second
term.\eject

It remains to deal with the second term of \eqref{separateone}.
Let $\delta>0$ be fixed (at first) and $m\dvtx \bbr\to]1,1+\delta[$ a
strictly monotone increasing auxiliary function. Since $m\geq1$ and
since we have at least one jump of size $>R$ on $\{\tau_R\leq t\}$, we obtain
\begin{eqnarray*}
\bbp^x(\tau_R\leq t) & \leq&\bbp^x \biggl(
\int_0^t \int_{\llVert
z \rrVert \geq R} m\bigl(
\llVert z \rrVert \bigr) \mu^{X^\sigma}(\cdot ;\mathrm{d}s,\mathrm{d}z) \geq m(R) \biggr)
\\
& \leq&\frac{1}{m(R)} \bbe^x \biggl( \int_0^t
\int_{\llVert  z
\rrVert \geq
R} m\bigl(\llVert z \rrVert \bigr)
1_{ [ [0, \sigma]  ]}(s) \mu^X(\cdot;\mathrm{d}s,\mathrm{d}z) \biggr)
\\
& =& \frac{1}{m(R)} \bbe^x \biggl( \int_0^t
\int_{z\neq0} m\bigl(\llVert z \rrVert \bigr)
1_{ [ [0 , \sigma[  [}(s) 1_{B_R(0)^c} (z) N(X_s,\mathrm{d}z) \, \mathrm{d}s \biggr)
\\
&\leq&(1+\delta) t \sup_{s\leq t\wedge\sigma} N\bigl(X_s,
B_R(0)^c\bigr)
\\
&\leq&(1+\delta) t \sup_{\llVert  y-x \rrVert \leq S} N\bigl(y,B_R(0)^c
\bigr)
\\
&\leq&(1+\delta) 4t \sup_{\llVert  y-x \rrVert \leq S} \int_{z\neq
0} \biggl(
\biggl\llVert \frac{z}{2R} \biggr\rrVert^2 \wedge1 \biggr)
N(y,\mathrm{d}z)
\end{eqnarray*}
because $m(\llVert  z \rrVert ) 1_{ [ [0 , \sigma[  [}(s)
1_{B_R(0)^c} (z)$ is in
class $F_p^1$ of Ikeda and Watanabe (see \cite{ikedawat}, Section
II.3). Since $\delta$ can be chosen arbitrarily small, we obtain by
Lemma \ref{lemoutsourced}
%
\begin{equation}
\label{finalterm2} \bbp^x(\tau_R\leq t) \leq4t
\sup_{\llVert  y-x \rrVert \leq S} \widetilde {c}_d \sup_{\llVert  \varepsilon \rrVert  \leq1} \biggl
\llvert p \biggl(y,\frac {\varepsilon}{2R} \biggr) \biggr\rrvert .
\end{equation}

Plugging together \eqref{finalterm1D}, \eqref{finalterm1DcA}, \eqref
{finalterm1DcB} and \eqref{finalterm2}, we obtain \eqref{techmainstopped}.

For the particular case $\sigma=\sigma_{3\widetilde{R}}^x$, we have
\[
\bigl\{ \bigl(X_\cdot^\sigma-x\bigr)_t^* \geq2
\widetilde{R} \bigr\} = \bigl\{ (X_\cdot -x)_t^* \geq2
\widetilde{R} \bigr\}
\]
and therefore, for every $\widetilde{R}>0$
%
\begin{equation}
\label{techmain} \bbp^x\bigl((X_\cdot-x)_t^*
\geq2\widetilde{R}\bigr) \leq c_d \cdot t \cdot \sup_{\llVert  y-x \rrVert \leq3\widetilde{R}}
\sup_{\llVert
\varepsilon \rrVert \leq
1} \biggl\llvert p \biggl(y,\frac{\varepsilon}{2\widetilde{R}} \biggr) \biggr
\rrvert .
\end{equation}
Setting $R:=(1/2)\widetilde{R}$, we obtain \eqref{firstestimate}. The
proof of \eqref{secondestimate} works literally as in the case of rich
Feller processes satisfying \eqref{growth} and \eqref{sector}.
Compare in this context \cite{schilling98}, Lemma 6.3 and Lemma 4.1.
The condition \eqref{growth} is not used in the proofs of these lemmas.
\end{pf*}

\begin{pf*}{Proof of Theorems \ref{thmttoinfty} and \ref{thmtto0}}
Since the proofs of the analogue statements for rich Feller processes
can be adapted and since all eight proofs are very similar, we decided
to give only on exemplary proof, namely of \eqref{proved}: Fix $x\in
\bbr^d$. Let $\lambda>\beta_\infty^x$ and choose $\lambda>\alpha_1>\alpha_2>\beta_\infty^x$. We have
\[
\bbp^x \bigl((X_\cdot- x)_t^* \geq
t^{1/\alpha_1} \bigr) \leq c_d \cdot t \cdot H\bigl(x,t^{1/\alpha_1}
\bigr) \leq c_d' \cdot t \bigl(t^{1/\alpha_1}
\bigr)^{-\alpha_2} =c_d' t^{1-(\alpha_2/\alpha_1)}
\]
for $t$ small enough, say $t<T_0$, since the $\lim\sup$ is considered.
Now let $t_k:=(1/2)^k$ for $k\in\bbn$. We obtain
\[
\sum_{k=k_0}^\infty\bbp^x \bigl(
(X_\cdot-x)_{t_k}^* \geq t_k^{1/\alpha_1}
\bigr) \leq c_d' \sum_{k=k_0}^\infty2^{-k(1-(\alpha
_2/\alpha_1))}
<\infty,
\]
where $k_0$ depends on $T_0$. By the Borel--Cantelli lemma, we obtain
\[
\bbp^x \Bigl( \mathop{\lim\sup}_{k\to\infty} (X_\cdot-x)_{t_k}^*
\geq (t_k)^{1/\alpha_1} \Bigr) =0
\]
and hence $(X_\cdot-x)_{t_k}^* < (t_k)^{1/\alpha_1}$ for all $k\geq
k_1(\omega)$ on a set of probability one. For fixed $\omega$ in this
set and $t_{k+1}\leq t \leq t_k$ and $k\geq k_1(\omega) \geq k_0$, we have
\[
\bigl(X_\cdot(\omega)-x\bigr)_t^* \leq
\bigl(X_\cdot(\omega)-x\bigr)_{t_k}^* \leq
t_k^{1/\alpha_1} \leq2^{1/\alpha_1} t^{1/\alpha_1}
\]
and since $\lambda>\alpha_1$
\[
t^{-1/\lambda} (X_\cdot-x)_t^* \leq2^{1/\alpha_1}
t^{(1/\alpha
_1)-(1/\lambda)}
\]
which converges $\bbp^x$-a.s to zero for $t\downarrow0$.
\end{pf*}

\section*{Acknowledgements}

The author wishes to thank Ren\'e Schilling (TU Dresden) for suggesting
the problem and Peter Furlan (TU Dortmund) for interesting discussions
on ODEs. Furthermore, he wishes to thank an anonymous referee for
his/her work. Financial support by the German Science Foundation (DFG)
for the project SCHN1231/1-1 is gratefully acknowledged.



\printhistory

\end{document}